\documentclass[11pt]{amsart}
\usepackage{geometry}             
\usepackage{graphicx}
\usepackage{amssymb}
\usepackage{epstopdf}
\usepackage{epstopdf}
\usepackage{psfrag}
\newtheorem{lemma}{Lemma}[section]
\newtheorem{theorem}{Theorem}[section]

\newtheorem{prop}{Proposition}[section]
\newtheorem{remark}{Remark}[section]

\renewcommand{\v}{ {\bf w}}

\newcommand{\W}{ {\bf W}}
\newcommand{\h}{ {\bf h}}
\newcommand{\f}{ {\bf f}}
\newcommand{\F}{ {\bf F}}
\newcommand{\oV}{ {\mathcal{V}}_i}
\newcommand{\be}{ {\mathbb{Z}} }
\newcommand{\done}{ {\bf d}_1 }
\begin{document}
\author{Pierre Degond}
\address[P. Degond]{Department of Mathematics,
Imperial College London,
London SW7 2AZ, United Kingdom}
\email{pdegond@imperial.ac.uk}
\author{Michael Herty}
\address[M. Herty]{RWTH Aachen University, Department of Mathematics,  
RWTH Aachen University,  D-52062 Aachen,  Germany ({\tt herty@igpm.rwth-aachen.de})}
\email{herty@igpm.rwth-aachen.de}
\author{Jian-Guo Liu}
\address[J.-G. Liu]{Department of Physics and Department of Mathematics, Duke University, Durham, NC 27708, USA
({\tt jliu@phy.duke.edu})}
\email{jliu@phy.duke.edu}
\title{  Meanfield games and model predictive control   }
\date{\today}       

\begin{abstract}
Mean-Field Games are games with a continuum of players that incorporate the time-dimension through a control-theoretic approach. Recently, simpler approaches relying on the Best Reply Strategy have been proposed. They assume that the agents navigate their strategies towards their goal by taking the direction of steepest descent of their cost function (i.e. the opposite of the utility function). In this paper, we explore the link between Mean-Field Games and the Best Reply Strategy approach. This is done by introducing a Model Predictive Control framework, which consists of setting the Mean-Field Game over a short time interval which recedes as time moves on. We show that the Model Predictive Control offers a compromise between a possibly unrealistic Mean-Field Game approach and the sub-optimal Best Reply Strategy. 
\end{abstract}
\maketitle
\section{Introduction} \label{sec:intro}

According to the definition of the Handbook \cite{FouqueLangsam2013aa}, systemic risk is the risk of a disruption of the proper functioning of the market which results in the reduction of the growth of the world's Gross Domestic Product (GDP). In economics, a system such as a market can be described by a game, i.e.  a set of agents endowed with strategies (and possibly other attributes) that they may play upon to maximize their utility function. In a game, the utility function depends on the other agents' strategies. The proper functioning of a market is associated to a Nash equilibrium of this game, i.e. a set of strategies such that no agent can improve on his utility function by changing his own strategy, given that the other agents' strategies are fixed. At the market scale, the number and diversity of agents is huge and it is more effective to use games with a continuum of players. Games with a continuum of players have been widely explored  \cite{Aumann1966aa,Mas-Colell1984aa,Schmeidler1973aa,ShapiroShapley1978aa}. 

To study systemic risk and its induced catastrophic changes in the economy, it is of primary importance to incorporate the time-dimension into the description of the system. A possible framework to achieve this is by means of a control-theoretic approach, where the optimal goal is not a simple Nash equilibrium, but a whole set of optimal trajectories of the agents in the strategy space. Such an approach has been formalized in the seminal work of \cite{LasryLions2007aa} and popularized under the name of 'Mean-Field Game (MFG)'. It has given rise to an abundant literature, among which (to cite only a few) \cite{BlanchetCarlier2014aa,BensoussanFrehseYam2014aa,BensoussanFrehseYam2013aa,GueantLasryLions2011aa,Cardaliaguet2010aa}. The MFG approach offers a promising route to investigate systemic risk. For instance, in the recent work \cite{CarmonaFouqueSun2013aa}, the MFG framework has been proposed to model systemic risk associated with inter-bank borrowing and lending strategies. 

However, the fact that the agents are able to optimize their trajectory over a large time horizon in the spirit of physical particles subjected to the least action principle can be seen as a bit unrealistic. A related but different approach has been proposed in \cite{DegondLiuRinghofer2012aa} and builds on earlier work on pedestrian dynamics \cite{DegondAppert-RollandMoussaid2013aa}. It consists in assuming that agents perform the so-called 'Best-Reply Strategy' (BRS): they determine a local (in time) direction of steepest descent (of the cost function, i.e. minus the utility function) and evolve their strategy variable in this direction. This approach has been applied to models describing the evolution of the wealth distribution among interacting economies, in the case of conservative \cite{DegondLiuRinghofer2014ac} and nonconservative economies \cite{DegondLiuRinghofer2014aa}. However, the link between MFG and BRS was still to be elaborated. This is the object of the present paper. We show that the BRS can be obtained as a MFG over a short interval of time which recedes as times evolves. This type of control is known as Model Predictive Control (MPC) or as Receding Horizon Control. The fact that the agents may be able to optimize the trajectories in the strategy space over a small but finite interval of time is certainly a reasonable assumption and this MPC strategy could be viewed as an improvement over the BRS and some kind of compromise between the BRS and a fully optimal but fairly unrealistic MFG strategy. We believe that MPC can lead to a promising route to model systemic risk. In this paper though, we propose a general framework to connect BRS to MFG through MPC and defer its application to specific models of systemic risk to future work. 

Recently, many contributions on meanfield games and control mechanisms for particle systems have been made. For more details on meanfield games we refer to  \cite{BlanchetCarlier2014aa,BensoussanFrehseYam2014aa,BensoussanFrehseYam2013aa,GueantLasryLions2011aa,Cardaliaguet2010aa,LasryLions2007aa}. Among the 
many possible meanfield games to consider we are interested in differential (Nash) games of possibly infinitely many particles (also called players). 
Most of the literature in this respect  
 treats theoretical and numerical approaches for solving the Hamilton--Jacobi Bellmann (HJB) equation for the value
function  of the underlying game, see e.g.  \cite{Cardaliaguet2010aa} for an overview. 
Solving the HJB equation allows to determine the optimal control for the particle game. However, the associated
 HJB equation posses several theoretical and numerical difficulties among which the need to solve it backwards in time is the most severe one, at least from 
 a numerical perspective. Therefore, recently model predictive control (MPC) concepts on the level particles or of the associated kinetic equation  have been proposed \cite{DegondAppert-RollandMoussaid2013aa,CouzinKrauseFranks2005aa,FornasierSolombrino2013aa,CaponigroFornasierPiccoli2013aa,AlbiHertyPareschi2014aa,DegondLiuRinghofer2014aa,DegondLiuRinghofer2014ac,CaponigroFornasierPiccoli2013aa}.
While MPC has been well established in the case of finite--dimensional problems \cite{GrunePannek2011ab,Sontag1998aa,MayneMichalska1990aa},  and also in engineering literature under the term receding horizon control, contributions to systems of infinitely many interacting particles and/or game theoretic questions related to infinitely many particles are rather recent. It has been shown that MPC concepts applied to problems of  infinitely many interacting particles have the advantage  to allow for efficient computation \cite{AlbiHertyPareschi2014aa,DegondLiuRinghofer2014ac}. However, by construction MPC only leads to suboptimal solutions, see for example \cite{HertySteffensenPareschi2014aa} for a comparison in the case of simple opinion formation model. Also, the existing approaches mostly for alignment models  do not necessarily treat game theoretic concepts but focus on 
for example sparse global controls \cite{CouzinKrauseFranks2005aa,FornasierSolombrino2013aa,CaponigroFornasierPiccoli2013aa}, time--scale separation and local mean--field controls 
\cite{DegondLiuRinghofer2014aa} called best--reply strategy, or MPC on very short time--scales \cite{AlbiHertyPareschi2014aa} called instantaneous control. 
Typically the MPC  strategy is obtained solving an auxiliary problem 
(implicit or explicit) and the resulting expression for the control is substituted back into 
the original dynamics leading to a possibly modified and new dynamics. Then, a meanfield description is derived using  Boltzmann or a macroscopic approximation. 
This requires the action of the control to be {\em local } in time  and independent of future states of the system contrary to solutions of the HJB equation. 
Usually in MPC approaches independent optimal control problems are solved where particles do not anticipate the optimal control choices
 other particles contrary to meanfield games \cite{LasryLions2007aa}.
\par
\begin{figure}[htb]\center
\includegraphics[width=.75\textwidth]{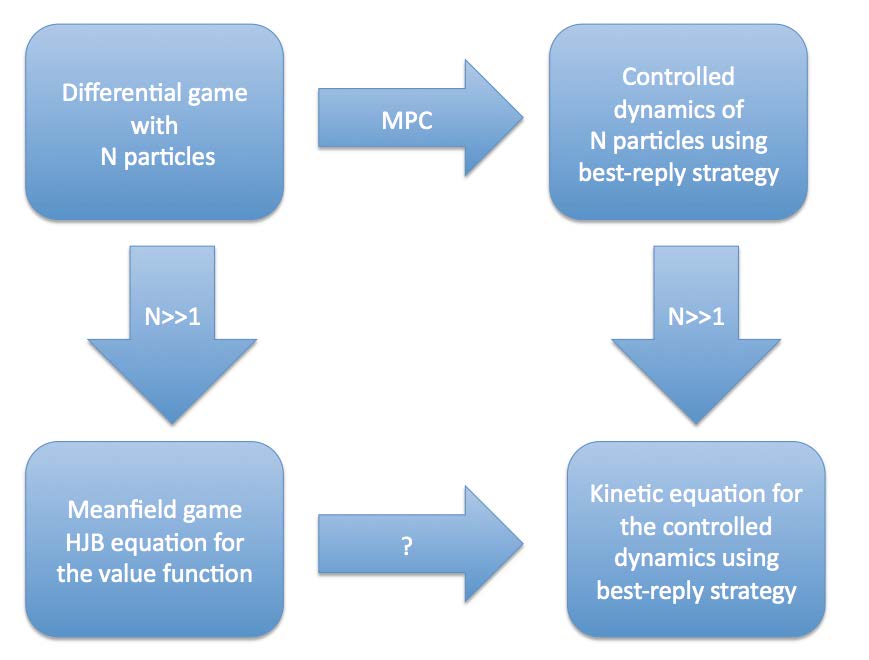}
\caption{ Relation between MPC concepts and meanfield games. The starting point are finite--dimensional differential games with $N$ players in the top left part (Section \ref{sec:setting}).
The connection for $N\to\infty$ of this games has been investigated for example in \cite{LasryLions2007aa,Cardaliaguet2010aa} 
and leads to the HJB for meanfield games in the bottom left part of the figure (Section \ref{sec:meanfield}). If applying MPC concepts to the differential game 
as for example the best-reply  strategy we obtain  a controlled dynamics for $N$ particles in the top right part \cite{DegondLiuRinghofer2014ac} (Section \ref{top-left-to-top-right}).
  The meanfield limit for $N\to \infty$ 
leads to a kinetic equation in the bottom right part (Section \ref{top-right-to-bottom-right}). Those results are summarized in Lemma \ref{lemma1}.
This paper also investigates the link between the meanfield game and the kinetic equation indicated 
by a question mark. The result is summarized in Proposition \ref{lemma2}.
 }
\label{fig1}
\end{figure}

 
In this paper we contribute to the recent discussion by formal computations leading to a link between meanfield games and MPC concepts
 proposed on the level of particle games and associated kinetic equations. The relationship we plan to establish is highlighted
in Figure \ref{fig1}. More precisely, we want to show that  the MPC concept of the best--reply strategy \cite{DegondLiuRinghofer2014ac} 
may be at least formally be derived from a meanfield games context. 
\section{Setting of the problem }\label{sec:setting}
\newcommand{\R}{\mathbb{R}}
We consider $N$ particles labeled by $i=1,\dots,N$ where each particle has a state 
$x_i \in \R.$ We denote by $X=(x_i)_{i=1}^N$ the state of all particles 
and by $X_{-i}=(x_j)_{j=1, j\not =i }^N$ the states of all particles except $i.$ Further,
we assume that each particle's dynamics is governed by  a smooth function $f_i:\R^N\to \R$
 depending on the state $X$ and we assume  that each particle 
may control its dynamics by a control $u_i.$ The dynamics for the particles $i=1,\dots,N$ 
is then given by 
\begin{equation}
\label{eq:full dynamics}
\frac{d}{dt} x_i(t) = f_i(X(t)) + u_i(t), \; i=1,\dots,N,
\end{equation}
and initial conditions 
\begin{equation}\label{eq:IC} 
x_i(0)=\bar{x_i}. 
\end{equation}
 We will drop the time--dependence of the variables whenever the intention is clear.
Examples of models of the type \eqref{eq:full dynamics} 
are alignment models in socio--ecological context, microscopic traffic flow models, production and many more, see e.g.
 the recent 
survey  \cite{MotschTadmor2013aa,NaldiPareschiToscani2010aa,VicsekZafeiris2012aa}. 
 In recent contributions to control theory for equation \eqref{eq:full dynamics}
the case of a {\em single} control variable $u_i \equiv u$ for all $i$ has been considered \cite{AlbiHertyPareschi2014aa,AlbiPareschiZanella2014aa,CaponigroFornasierPiccoli2013aa}. 
Here, we allow each particle to chose its own control strategy $u_i$. We suppose a control horizon of $T>0$ be given. 
As in \cite{Cardaliaguet2010aa} we suppose that particle $i$ minimizes its own objective functional 
 and determines therefore the optimal $u_i^{*}$ by 
\begin{equation}
\label{eq:general min particle}
u^{*}_i(\cdot) = \mbox{ argmin }_{ u_i: [0,T] \to \R } \int_0^T \left( \frac{ \alpha_i(s) }2  u_i^2(s) + h_i(X(s)) \right) ds,  \; i=1,\dots,N.
\end{equation}
Herein, $X(s)$ is the solution to \eqref{eq:full dynamics} and equation \eqref{eq:IC}. The optimal control and the corresponding 
optimal trajectory will from now on be denoted with  superscript $\ast.$
The minimization is performed on all sufficiently smooth functions $u_i:[0,T] \to \R$. There is no restriction on the control $u_i$
 similar to \cite{LasryLions2007aa}.  The objective  $h_i:\R^{N}\to \R$ related
to particle $i$ is also supposed to be sufficiently smooth. The weights of the control $\alpha_i(t) >0, \forall i, \; t\geq 0$ and  under additional conditions convexity
of each  optimization problem \eqref{eq:general min particle} is guaranteed.  As seen in Section \ref{sec:meanfield} a challenge in solving the problem \eqref{eq:general min particle} 
relies on the fact that the associated HJB has to be solved backwards in time.  Contrary to \cite{AlbiHertyPareschi2014aa,CaponigroFornasierPiccoli2013aa} 
problem \eqref{eq:general min particle} are in fact $N$ optimization problems that need to be solved {\em simultaneously } due to the dependence 
of $X$ on $U=(u_i)_{i=1}^{N}$ through equation \eqref{eq:full dynamics}.  This implies that each particle $i$ {\em anticipates the optimal } strategy
of all other particles $U^{*}_{-i}$ when determining its optimal control $u^{*}_i$. Obviously, the problem \eqref{eq:general min particle} simplifies
when each particle $i$ {\em anticipates an a priori fixed strategy} of all other particles $U_{-i}.$ Then, the problem \eqref{eq:general min particle}
decouples (in $i$) and the optimal strategy $u_i$ is determined independent of the optimal strategies $U^{*}_{-i}.$  It has been argued that 
this is the case for reaction in pedestrian motions \cite{DegondLiuRinghofer2014aa}.
In fact, therein the following best--reply strategy has been proposed as a substitute for problem \eqref{eq:full dynamics}
\begin{equation}\label{eq:best reply}
u_i(t) = - \partial_{x_i} h_i(X(t)), \; t \in [0,T].
\end{equation}
As in the meanfield theory presented in \cite{LasryLions2007aa,Cardaliaguet2010aa} we need to impose  assumptions {\bf (A)} on $f_i(X)$
and $h_i(X)$ before passing to the limit $N\to\infty.$ The assumption {\bf (B)} will be used in Section \ref{sec:meanfield}.  
\begin{itemize}
\item[{\bf (A)}] For all $i=1,\dots,N$ and any permutation $\sigma_i : \{1,\dots,N\} \backslash \{ i \} \to \{1,\dots,N\} \backslash \{ i \}  $
we have 
$$f_i(X) = f(x_i,X_{-i})  \mbox{ and  }  f(x_i,X_{-i}) = f(x_i, X_{\sigma_i} )$$ 
for a smooth function $f:\R \times \R^{N-1} \to \R$ and where $X_{\sigma_i} = (x_{\sigma_i(j)})_{j=1, j\not = i}^{N}.$
 Further we assume that for each $i$ the function $h_i(X)$ enjoys  the same properties as stated for $f_i(X).$ 
 \item[ {\bf (B)}] We assume that $\alpha_i(t) = \alpha(t)$ for all $t \in [0,T]$ and all $i=1,\dots,N.$  
\end{itemize}
Under additional growth conditions sequences of symmetric functions in many variables
have a  limit in the space of functions defined on probability measures, see e.g. 
\cite[Theorem 2.1]{Cardaliaguet2010aa}, \cite[Theorem 4.1]{BlanchetCarlier2014aa}. 
The corresponding result is recalled as Theorem~\ref{Theorem2.1Card} in the appendix 
for convenience. 
\par 
To exemplify computations later
on we will use a basic wealth model
 \cite{BouchaudMezard2000aa,CordierPareschiToscani2005aa,DegondLiuRinghofer2012aa,DegondLiuRinghofer2014aa}
where 
\begin{equation}
\label{eq:ex:opinion}
f_i(X) = \frac{1}N \sum\limits_{j=1}^N P(x_i,x_j) ( x_j - x_i)
\end{equation}
for some bounded, non--negative and smooth function $P(x,\tilde{x}).$ Clearly, $f$ fulfills {\bf(A)}.
 As objective function we use a measure depending only on aggregated quantities as in \cite{DegondLiuRinghofer2014ac}. 
 An example fulfilling {\bf (A)} is 
 \begin{equation}
 \label{eq:ex:objective}
 h_i(X) = \frac{1}{N-1} \sum\limits_{j=1, j\not = i}^{N} \phi(x_i,x_j)
 \end{equation}
for some smooth function $\phi:\R\times \R \to\R.$ 
\par
Finally, we introduce some additional notation. We denote by $\mathcal{P}(\R)$ the space of Borel probability measures over $\R$. 
The empirical discrete probability measure $m^{N} \in \mathcal{P}(\R)$ concentrated at a positions $X \in \R^N$ is denoted by 
$$ m^{N}_X = \frac{1}N \sum\limits_{i=1}^N \delta(x-x_i).$$
We also use this notation if $X$ is time dependent, i.e., $X=X(t)$, leading to the family of probability measures $m^{N}_X = m^{N}_X(t)=\frac{1}N \sum\limits_{i=1}^N \delta( x- x_i(t)).$ 
If the intention is clear we do not explicitly denote the dependence on $x$ of the measure $m^{N}_X$ (and on time $t$ if $X=X(t)$ is time-dependent).
\par 
Based on the assumption {\bf (A)} we will frequently use Theorem \cite[Theorem 2.1]{Cardaliaguet2010aa}. 
This theorem is repeated for convenience in the appendix as Theorem \ref{Theorem2.1Card}: let  $g:\R^{N}\to \R$ such that $g$ is symmetric $g(X)=g(X_\sigma)$
where $X_\sigma=(x_{\sigma(i)})_{i=1}^{N}$ and any permutation $\sigma: \{ 1,\dots, N\} \to \{1,\dots,N\}$. We may extend $g$ to a function $g^{N}:\mathcal{P}(\R) \to \R$ such that  
$g^{N}(m^{N}_X) = g(X).$ Under  assumptions given in Theorem \ref{Theorem2.1Card} the 
family $(g^{N})_{N=1}^\infty$ is equicontinuous and there exists a limit ${\bf g}:\mathcal{P}(\R)\to \R$ such that up to a subsequence
$\lim\limits_{N\to \infty} \sup_{X} | g(X) - {\bf g}( m^{N}_X ) | = 0.$
The result can be extended to a family of functions $f$ fulfilling assumption {\bf (A)}, see \cite[Theorem 4.1]{BlanchetCarlier2014aa}.  
We obtain an equicontinuous family $(f^{N})_{N=1}^\infty$ with $f^{N}:\R \times \mathcal{P}(\R) \to \R$ such that 
$f^{N}(\xi, m^{N-1}_{X_{-i}}) = f(\xi, X_{-i})$ with limit $\f:\R \times \mathcal{P}(\R) \to \R$ and such that for any $i \in \{ 1,\dots,N \}$ and a compact set $Q \subset \R^{N-1}$ 
we have for any fixed $R>0$
\begin{equation}
\label{def:con} \lim\limits_{N\to \infty} \sup_{ | \xi | < R, X_{-i} \subset Q } | f(\xi,X_{-i}) -  \f(\xi,m^{N}_X) | =  0.
\end{equation}
In equation \eqref{def:con}  we have the empirical measure on $N$ points in the argument of $\f$ even so $f^{N}$ is defined on the empirical measure $m^{N-1}_{X_{-i}}$, i.e., 
we have  $$\lim\limits_{N\to \infty} \sup_{ | \xi | < R, X_{-i} \subset Q } | f^{N}(\xi,m^{N-1}_{X_{-i}}) -  \f(\xi,m^{N}_X) | =0.$$ 
 This is true since in the definition of $f^{N}$  the contribution of the empirical measure is $\frac{1}N$ for each point $x_i.$ More details are given in \cite{BlanchetCarlier2014aa}
and \cite[Section 7]{Cardaliaguet2010aa}. Since in the following it is often of importance to highlight the dependence on the empirical measure $m^{N}_X$
we introduce the following notation: we write 
$$ f(\xi,X_{-i}) = f^{N}(\xi, m^{N-1}_{X_{-i}} ) \sim \f(\xi,m^{N}_X), \; $$
whenever equation \eqref{def:con} holds true. 

\subsection{From differential games to controlled particle dynamics } \label{top-left-to-top-right}
The best--reply strategy \eqref{eq:best reply} is obtained also from  a MPC approach \cite{MayneMichalska1990aa} 
applied to equations \eqref{eq:full dynamics} and \eqref{eq:general min particle}. In order to derive the best--reply 
strategy we consider the following problem: suppose we are given the state $X(t)$ of the system \eqref{eq:full dynamics} 
at time $t>0$. Then, we consider a control horizon of the MPC of $\Delta t>0$ and supposedly small. We assume
that the applied control $u_i(s)$ on $(t, t+\Delta t)$ is constant. For particle $i$ we 
denote the unknown constant by $\tilde{u}_i.$  Instead of solving the  problem 
\eqref{eq:general min particle} now on the full time interval $(t,T)$ 
we consider the objective function only on the receding time horizon $(t,t+\Delta t).$ 
Further, we discretize the dynamics \eqref{eq:full dynamics} on $(t,t+\Delta t)$
 using an explicit Euler discretization for the initial value $\bar{X} = X_i(t).$ We discretize the objective function
 by a Riemann sum. A naive discretization leads to a penalization of the control of the type 
 $\frac{\alpha_i(t+\Delta t)}2 \tilde{u}^{2}.$ Since the explicit Euler discretization in equation \eqref{eq:MPC 1} is only accurate up to order $O( (\Delta t)^{2} )$ 
we additionally require to have $\tilde{u}_i = O(1)$ to obtain a meaningful control in the discretization \eqref{eq:MPC 1} and also in the limit for $\Delta t \to 0$. 
Therefore, in the MPC problem we need to scale the penalization of the control accordingly by $\Delta t.$  
 This leads to a MPC problem associated with equation \eqref{eq:general min particle} and  given by 
 \begin{align}
 \label{eq:MPC 1} x_i(t+\Delta t) = \bar{x}_i +\Delta t  \left( f_i( \bar{X} ) +  \tilde{u}_i \right), \;& & i=1,\dots, N, \\
 \label{eq:MPC 2} \tilde{u}_i = \mbox{ argmin }_{ \tilde{u} \in \R}  \Delta t \left(  h_i\left( X(t+\Delta t) \right)  + \Delta t  \frac{ \alpha_i( t+ \Delta t) }2 \tilde{u}^2 \right)  , \;& &  i=1,\dots, N.
 \end{align}
 Solving the minimization problem \eqref{eq:MPC 2} leads to 
\begin{equation*}
\alpha_i(t+\Delta t) \; \tilde{u}_i = - \partial_{x_i} h_i( \bar{X} ), \; i=1,\dots,N.
\end{equation*}
Now, we obtain a $\tilde{u}_i$ of order $O(1)$ by Taylor expansion of $\alpha_i$ at time $t.$ 
Within the MPC approach the control for the time interval $(t,t+\Delta t)$ is therefore given by   equation \eqref{eq:MPC best reply}. 
\begin{equation}\label{eq:MPC best reply}
 \tilde{u}_i = - \frac{1}{\alpha_i(t)} \partial_{x_i} h_i( \bar{X} ), \; i=1,\dots,N.
\end{equation}
Usually, the dynamics \eqref{eq:MPC 1} is then computed with the computed  control up to $t+\Delta t.$ Then, the process is repeated
using the new state $X(t+\Delta t).$ Substituting \eqref{eq:MPC best reply} into \eqref{eq:MPC 1} and letting $\Delta t \to 0$ 
we obtain 
\begin{equation}\label{eq:MPC controlled dynamics}
\frac{d}{dt} x_i(t) = f_i(X(t))  - \frac{1}{\alpha_i(t)} \partial_{x_i} h_i( X(t) ), \; i=1,\dots, N, t \in [0,T].
\end{equation}
This dynamics coincide with the dynamics generated by the best--reply strategy \eqref{eq:best reply} provided that $\alpha_i(t) \equiv 1.$ 
Therefore, on  a particle level the controlled dynamics \eqref{eq:MPC controlled dynamics} of the best--reply strategy \cite{DegondLiuRinghofer2014ac} is equivalent to a MPC 
formulation of the problem \eqref{eq:general min particle}. For the toy example we obtain 
\begin{equation}
\label{eq:ex:MPC} \frac{d}{dt} x_i =  \frac{1}N \sum\limits_{j=1}^N P(x_i,x_j) ( x_j - x_i) - \frac{1}{ (N-1) \alpha_i(t)}  \sum \limits_{j=1, j\not = i}^N \partial_{x_i} \phi(x_i, x_j).
\end{equation}

\subsection{From controlled particle dynamics \eqref{eq:MPC controlled dynamics} to kinetic equation} \label{top-right-to-bottom-right}

The considerations herein have essentially been studied for the best--reply strategy in the series of papers \cite{DegondLiuRinghofer2012aa,DegondLiuRinghofer2014ac,DegondLiuRinghofer2014aa} and it is only repeated for 
convenience. The starting point is the controlled dynamics given by equation \eqref{eq:MPC controlled dynamics} which slightly extends
 the best--reply 
strategy. In order to pass to the meanfield limit we assume that {\bf (A)} and {\bf (B)} holds true. Then the particles are governed by 
\begin{equation} 
\label{eq:controlled dynamics}
\frac{d}{dt} x_i(t) = f(x_i(t), X_{-i}(t)) - \frac{1}{\alpha(t)} \partial_{x_i} h(x_i(t), X_{-i}(t) ), \; i=1,\dots, N.
\end{equation}

Associated with the trajectories $X=X(t)$ the discrete probability measure $m^N_X$ is given by $m^{N}_X=\frac{1}N \sum\limits_{j=1}^N \delta(x-x_i(t)).$ 
Using the weak formulation for a test function $\psi:\R \to \R$ we compute the dynamics of $m^N_X$ over time as 
\begin{align*}
\frac{d}{dt} \int \psi(x) m^N_X dx = \frac{1}N \sum\limits_{i=1}^N \int \psi'(x)  \left( f(x,X_{-i}) - \frac{1}{\alpha} \partial_x h(x,X_{-i}) \right) \delta(x-x_i(t)) dx
\end{align*}

Using \cite[Theorem 2.1]{Cardaliaguet2010aa}  and denoting by $m^{N-1}_{X_{-j}}(t) = \frac{1}{ N-1} \sum\limits_{k=1, k\not =j } \delta(x-x_k(t))$ 
a family of  empirical measures on $\R$  
 we obtain from the previous equation 
\begin{align*}
\frac{d}{dt} \int \psi(x) m^N_X dx = \frac{1}N \sum\limits_{i=1}^N \int \psi'(x)  \left( f^N(x, m^{N-1}_{X_{-i}} ) - 
\frac{1}{\alpha} \partial_x h^N(x,m^{N-1}_{X_{-i}} ) \right)  \delta(x-x_i(t)) dx
\end{align*}
for  some function $f^N, h^N:\R \times \mathcal{P}(\R) \to \R$.  Assume $f$ and $h$ 
fulfill the assertions of \cite[Theorem 4.1]{BlanchetCarlier2014aa}.    Then,   we obtain in 
the sense of equation \eqref{def:con} that for any $i$  and $N$ sufficiently large
$\f(x,m^{N}_X) \sim f^{N}(x,m^{N-1}_{X_{-i}})$ and  $\h(x,m^{N}_X) \sim h^{N}(x,X_{-i})$ and therefore
in the sense of equation \eqref{def:con}
\begin{align*}
\frac{d}{dt} \int \psi(x) m^N_X dx = \frac{1}N \sum\limits_{i=1}^N \int \psi'(x)\left( {\bf f}(x, m^N_X ) - 
\frac{1}{\alpha} \partial_x {\bf h}(x,m^N_X ) \right)   \delta(x-x_i(t))  dx =  \\ 
\int \psi'(x) m^N_X \left( {\bf f}(x, m^N_X ) -  \frac{1}{\alpha} \partial_x {\bf h}(x,m^N_X ) \right) dx  
\end{align*}
 
This is the weak form of the kinetic equation for a probability measure $m=m(t,x)$ 
\begin{equation}\label{eq:best reply kinetic}
\partial_t m + \partial_x \left( m \left(  {\bf f}(x,m) - \frac{1}\alpha \partial_x {\bf  h } (x,m) \right) \right) = 0
\end{equation}

\par
\noindent
For the toy example the corresponding function $f^N$ is given 
by $$f^N(x,m^N_{X_{-i}} ) = \frac{N-1}N \left( \int P(x,y) (y-x) m^N_{X_{-i}}(t,y)  dy \right) $$
and ${\bf f}(x,m^N_X) = \int P(x,y) (y-x) m^N_X(t,y)  dy.$ 
Therefore, we obtain  
\begin{equation*}
\partial_t m(t,x) + \partial_x \left(  \int  \left( P(x,y)(y-x)  - \frac{1}\alpha \partial_x \phi(x,y) \right) m(t,x) m(t,y) dy \right) = 0.
\end{equation*}


We summarize the previous findings in the following lemma.

\begin{lemma}\label{lemma1}
Consider a fixed time horizon $T>0$ and consider $N$ particles governed by the dynamics \eqref{eq:full dynamics} and initial conditions \eqref{eq:IC}. Assume {\bf (A)} and {\bf (B)}
hold true.  Assume each particle $i=1,\dots,N$, chooses its  control $u_i$ at time $t \in [0,T]$ by 
\begin{equation}\label{lem:ctrl} u_i(t) = - \frac{1}{\alpha_i(t)} \partial_{x_i} h_i(X(t)).\end{equation}
Assume that $m^{N}_{\bar{X}} \to {\overline{m}} \in \mathcal{P}(\R)$ for $N\to\infty.$ Then, 
the meanfield limit of the  particle dynamics \eqref{eq:full dynamics} and \eqref{lem:ctrl}  is given by 
\begin{equation}
\label{lem:mf} 
\partial_t m + \partial_x \left( m \left(  {\bf f}(x,m) - \frac{1}\alpha \partial_x {\bf  h } (x,m) \right) \right) = 0
\end{equation}
for initial data $m(t=0,x)=\overline{m}.$ 
\end{lemma}

Finally, we summarize the MPC approach. Consider a time interval $\Delta t >0$, an equidistant discretization $t_\ell = (\ell-1) \Delta t, \ell=1,\dots, N_T$  such that $N_T \Delta t =T$  and 
 a first--order numerical discretization of $N$ particle dynamics \eqref{eq:full dynamics} given by 
\begin{equation}\label{lem:disc}
x_{\ell,i} = x_{\ell-1,i} + \Delta t \left( f_i(X(t_{\ell})) + \tilde{u}_{\ell,i} \right), \; \ell = 1,\dots, N_T \mbox{ and }  \; x_{0,i}=\bar{x}_i.
\end{equation}
for $x_i(t_\ell)=x_{\ell,i} $ and $i=1,\dots, N.$ Let in equation \eqref{lem:disc} be the control $u_i(t) = \sum\limits_{\ell=1}^{N_T} \chi_{  \left[ (\ell-1)\Delta t, \ell \Delta t \right) }(t) \tilde{u}_{\ell,i}$ is piecewise constant. The constants $\tilde{u}_{\ell,i}$ obtained as discretization in time  of  equation \eqref{lem:ctrl} are given by 
\begin{equation}\label{lem:discctl} \tilde{u}_{\ell,i} = - \frac{1}{\alpha_i(t_\ell)} \partial_{x_i} h_i(X(t_\ell)), \; \ell=1,\dots,N_T.
\end{equation}
Then, each $\tilde{u}_{\ell,i}$ coincides up to $O(\Delta t)$ with optimal control on the time interval $\left[ (\ell-1)\Delta t, \ell \Delta t \right)$ of the 
 minimization problems \eqref{lem:min pb}. For each $\ell=1,\dots,N_T$ and each fixed value $X(t_{\ell-1})$ we have up $O(\Delta t)$ 
\begin{equation}
\label{lem:min pb}
\tilde{u}_{\ell,i}  = \mbox{argmin}_{u \in \R} \left(  h_i(X(t_\ell)) + \Delta t \; \frac{\alpha_i(t_\ell)}2 u^2 \right), i=1,\dots, N, 
\end{equation}
where $X(t_\ell)$ is given by equation \eqref{lem:disc}.

\section{Results related to meanfield games }\label{sec:meanfield}
In this paragraph we consider the limit of the problem \eqref{eq:general min particle}
for a large number of particles. This has been investigated for example in \cite{LasryLions2007aa}
and derivations (in a slightly different setting) have been detailed in \cite[Section 7]{Cardaliaguet2010aa}. In order to show the links presented in Figure \ref{fig1} 
we repeat the formal computations in \cite{Cardaliaguet2010aa}.
\par 
A notion of solution to the competing $N$ optimization problem \eqref{eq:general min particle} 
is the concept of Nash equilibrium. If it exists it may be computed for the differential games using the HJB equation. We briefly present computations leading to the HJB equation. Then, we discuss
the large particle limit of the HJB equation and derive the best--reply strategy.

\subsection{Derivation of the finite--dimensional HJB equation }
The HJB equation describes the evolution of a value function $V_i=V_i(t,Y)$ 
of the particle $i.$  The value function is defined as the future costs for a 
particle trajectory governed by equation \eqref{eq:full dynamics} and starting
at time $t \in (0,T)$, position $Y$ and  control $u_i:(t,T)\to \R, i=1,\dots,N,$ 
\begin{equation}
\label{eq:value fct i}
V_i(t, Y ) = \int_t^T \left( \frac{ \alpha_i(s)}2 u_i^2(s) + h_i(X(s))  \right)ds,
\end{equation}
where $X(s)=(x_i(s))_{i=1}^N$ is the solution to equation \eqref{eq:full dynamics} with control 
$U$ and initial condition
\begin{equation}
\label{eq:ic value i}
X(t) = Y.
\end{equation}
Among all possible controls $u_i$ we denote by $u_i^{*}$ the 
 optimal control that minimizes $V_i(t,Y)$. We investigate the relation of the value function of particle $i$ to 
the optimal control $u_i^*$. To this end assume that the coupled problem \eqref{eq:general min particle} 
has a unique solution denoted by $U^*=(u^*_i)_{i=1}^N$. 
Each  $u_i^* : [t,T]\to \R$ for each $i=1,\dots,N$, is hence a solution to 
$$ u^{*}_i = \mbox{ argmin }_{ u_i(\cdot):[t,T]\to\R } \{ V_i(t,Y): X \mbox{ solves } \eqref{3:state} \}, \; i=1,\dots,N. $$
The corresponding particle trajectories are denoted by $X^* = (x_i^*)_{i=1}^N$
and are obtained through \eqref{eq:full dynamics} for an initial condition  $X^*(0)=\bar{X}.$
\par 
Since $X^*(\cdot)$ depends on $U^*,$ minimizing the value function \eqref{eq:value fct i}
leads to the computation of formal derivatives of $V_i$ with respect to $u_i.$ The optimal 
control $u_i^*$ is then found as formal point (in function space) where the  derivative
of $V_i$ with respect to $u_i$ vanishes. We have as Gateaux derivative of $V_i$ in an arbitrary direction $v:[0,T]\to \R$ 
\begin{align}\label{3:reduced cost}
\frac{d}{d u_i } V_i(t,Y)[v] = \int_t^T \left( \alpha_i(s) u_i^*(s)  + \sum\limits_{k=1}^N \partial_{x_k} h_i(X^*(s)) \partial_{u_i} \left( x_k^{*}(s) \right)   \right) v(s) ds = 0 
\end{align}
The derivative is not easily computed due to the unknown derivative of each state $x_k^{*}$ with respect to the acting control $u_i^{*}.$ 
However, choosing a set of suitable co-states $\phi_j^\ell:[0,T] \to \R$ for $\ell=1,\dots,N$ and $j=1,\dots,N,$ we may simplify the 
previous equation \eqref{3:reduced cost}: we  test equation \eqref{eq:full dynamics} by  functions $\phi_j^\ell :[0,T] \to \R$ 
for $\ell,j=1,\dots,N$ such that 
$\phi_j^{\ell}(T)=0$ and integrate on $(t,T)$ with $0\leq t<T$ and  initial data at $X^*(t)=Y$ to obtain 
\begin{align}\label{3:weak state}
 \sum\limits_{j=1}^N \left\{  \int_t^T - \frac{d}{ds}  (\phi_j^\ell(s) ) x_j^*(s) - \phi_j^\ell(s)  \left( f_j(X^*(s)) + u_j^*(s) \right)  ds - \phi_j^\ell(t)y_j \right\} = 0, \; \ell=1,\dots,N.
\end{align} 
The derivative with respect to $u_i$ in an arbitrary direction $v$ is then 
\begin{align} \nonumber 
 \int_t^T \left\{  \sum\limits_{j=1}^N \left( - \frac{d}{ds}  (\phi_j^\ell(s) )  \partial_{u_i} \left( x_j^*(s) \right)  - \phi_j^\ell(s)  \left(  \sum\limits_{k=1}^N  \partial_{x_k} \left( f_j(X^*(s)) \right) \partial_{u_i} (x_k^*(s))  \right)  \right) - \phi_i^\ell(s)    \right\} \\ 
 \times  v(s)  ds = 0.\label{3:weak derivative}
\end{align} 
The previous equation can be equivalently rewritten as
\begin{align}  
 \sum\limits_{k=1}^N   \left(    - \frac{d}{ds} \phi^\ell_k(s) - \sum\limits_{j=1}^N \phi_j^\ell(s) \partial_{x_k} f_j(X^*(s)) \right) \partial_{u_i} ( x_k^*(s) ) v(s)  
ds =  \int_t^T \phi_i^\ell(s) v(s) ds.\label{eq:p-ast}
\end{align}
Let $\phi^i_j$ for $i,j=1,\dots,N$   fulfill the coupled linear system of adjoint equations  (or co-state equations), solved backwards in time,
\begin{align}\label{3:adjoint}
- \frac{d}{dt} \phi^i_j(t) - \sum\limits_{k=1}^N \phi^i_k (t) \partial_{x_j} ( f_k(X^*(t)) = \partial_{x_j} h_i(X^*(t)), \; \phi^i_j(T)=0. 
\end{align}
Then, formally for every $s \in (t,T)$ we have
\begin{align*}
 \sum\limits_{j=1}^N \partial_{x_j} h_i(X^*(s)) \partial_{u_i} (x_j(s) ) = 
\sum\limits_{k=1}^N 
- \frac{d}{dt} \phi^i_k(s) \partial_{u_i} (x^{*}_k(s))  - \sum\limits_{j=1}^N \sum\limits_{k=1}^N \phi^i_j (s) \partial_{x_k} ( f_j(X^*(s)) \partial_{u_i} (x^{*}_k(s)). 
\end{align*}
Since for all $v$ we have 
\begin{align*}
\sum\limits_{j=1}^N \int_t^T  \partial_{x_j} h_i(X^*(s)) \partial_{u_i} \left( x_j^{*}(s) \right) v(s) ds = \int_t^T \phi_i^{i}(s) v(s) ds
\end{align*}
it follows that $$\sum\limits_{j=1}^N \partial_{x_j} h_i(X^*(s)) \partial_{u_i} \left( x_j^{*}(s) \right) = \phi^{i}_i(s)$$ for $s \in (t,T).$ At
optimality the necessary condition is 
$$\frac{d}{du_i} V_i(t,Y)[v] = 0$$ 
for all $v$ which implies that thanks to equation \eqref{3:reduced cost} we obtain for a.e. $s\in (t,T)$
$$ 
\left( \alpha_i(s) u_i^*(s) \right) + \sum\limits_{j=1}^N  \partial_{x_j} h_i(X^*(s)) \partial_{u_i} \left( x_j^{*}(s) \right) =0. 
$$
This leads to the following equation a.e. in $s \in (t,T)$
\begin{align}
\label{3:gradient}
\alpha_i(s) u_i^*(s) + \phi_i^i(s) = 0.
\end{align} 
Equations 
\begin{equation}\label{3:state} 
\frac{d}{ds} x_j^*(s) = f_j(X^*(s)) + u_j^*(s), \; X^*(t)=Y.
\end{equation}
and  \eqref{3:adjoint} for $j=1,\dots,N$ and equation \eqref{3:gradient} comprise the optimality
conditions for the minimization of the value function $V_i(t,Y)$ given by equation \eqref{eq:value fct i}. Due to 
equation \eqref{3:state} this system is a coupled system of ordinary differential and algebraic equations
in the unknowns $\mathcal{S}:=(x_j^*, u_j^*, (\phi_j^i)_{i=1}^N )_{j=1}^N.$ Solving for all those unknowns yields
in particular the optimal control $u_i^*$ for the value function $V_i(t,Y)$  for all $i=1,\dots,N.$ The adjoint equation \eqref{3:adjoint}
is posed backwards in time such that the system is two--point boundary value problem and due to the strong coupling
of $x_j^*$ and $\phi^i_j$ it is not easy to solve. The derived system is a version of Pontryagins maximum principle  (PMP)
for sufficient regular and unique controls. We refer to \cite{Friedman1974aa,Sontag1998aa,Bressan2011aa} for more details.
From now on we assume that equation \eqref{3:gradient} where $\phi_j^i$ solves equation \eqref{3:adjoint} is
 necessary and sufficient for optimality of $u_i^*$ for minimizing the value function $V_i(t,Y).$ The corresponding
 optimal trajectory and co-state is denoted by $\mathcal{S}$ introduced above. We formally
 derive the HJB based on the previous equations of  PMP and refer to \cite[Chapter 8]{Friedman1974aa} 
 for a careful theoretical discussion.  
\par 
Consider the function $V_i(t,Y)$ evaluated along the optimal trajectory $\mathcal{S},$ i.e., 
let $\oV(t)=V_i(t, X^*(t)).$ Then, by definition of $V_i$ and $\mathcal{S}$ we have 
\begin{align*}
 - \frac{\alpha_i(t)}2  (u_i^*)^2(t) - h_i(X^*(t))  = \frac{d}{dt} \oV(t)  \\
=  \partial_t V_i(t,X^*(t)) + \sum\limits_{k=1}^{N} \partial_{x_k} V_i(t,X^*(t)) \left( f_k(X^*(t))+u_k^*(t) \right).
\end{align*}
Using the necessary condition \eqref{3:gradient} we obtain 
\begin{align}\label{3:temp}
 - \frac{1}{ 2\alpha_i(t)}  (\phi^i_i)^2(t) - h_i(X^*(t))  \\
 = \partial_t V_i(t,X^*(t)) + \sum\limits_{k=1}^{N} \partial_{x_k} V_i(t,X^*(t)) \left( f_k(X^*(t))-\frac{1}{\alpha_k(t)} \phi_k^k(t) \right). \nonumber
\end{align}
The trajectory of $X^*(s)$ depends on the initial condition $Y=(y_i)_{i=1}^N$. 
Computing the variation of $V_i(t,Y)$ with respect to $y_o$ for $o \in \{1,\dots,N\}$ and evaluating at $\mathcal{S}$ yields (since $u_i^{*}$ does not dependent
on $Y$):
\begin{align*}
\partial_{y_o} V_i(t,Y) = \int_t^T \left( \sum\limits_{k=1}^N \partial_{x_k} h_i(X^*(s)) \partial_{ y_o } (x_k^*(s))  \right) ds.
\end{align*}
From the weak form of the state equation \eqref{3:weak state} we obtain after differentiation with respect to the initial condition $y_o$ for  $\ell=1,\dots,N$ 
\begin{align*}
\sum\limits_{j=1}^N \int_t^T  \left( -\frac{d}{ds} \phi^\ell_j(s)  \partial_{y_o} (x_j^*(s) ) - \phi_j^\ell(s) \sum\limits_{k=1}^N  (\partial_{x_k}  f_j) (X^*(s) ) \partial_{y_o} ( x_k^*(s) ) \right) 
ds = \phi^\ell_o(t).
\end{align*}
Similarly  to the computations before we use the equation for $\phi^i_j$ given by equation \eqref{3:adjoint} 
to express 
\begin{align*}
\partial_{y_o} V_i(t,Y)  = \int_t^T \left( \sum\limits_{k=1}^N  (\partial_{x_k}  h_i) (X^*(s))  \partial_{y_o} ( x_k^*(s) ) \right) ds =  \\
\int_t^T \left( \sum\limits_{k=1}^N \left( 
 - \frac{d}{ds} \phi_k^i(s)  - \sum\limits_{j=1}^N \phi_j^i(s) (\partial_{x_k}  f_j) (X^*(s) ) \right) \partial_{y_o} ( x_k^*(s) ) \right) ds  = \phi_o^i(t).
\end{align*}
Therefore, $\nabla_Y V_i(t,Y) = (\phi^{i}_k)_{k=1}^{N}$ provided that $\phi^{i}_k$ is a solution to equation \eqref{3:adjoint}.
Now, along $\mathcal{S}$ we may express in equation \eqref{3:temp} the co--state by the derivative of $V_i$ with respect to $Y.$ Replacing $Y=X^{*}(t)$ 
we obtain 
\begin{align*}
- \frac{1}{2\alpha_i(t)} \left( \partial_{x_i} V_i(t,X^{*}(t)) \right)^{2} - h_i(X^*(t)) = \partial_t V_i(t,X^{*}(t)) +  \\
\sum\limits_{k=1}^{N} \partial_{x_k} V_i(t,X^{*}(t)) \left( f_k(X^{*}(t)) - \frac{1}{\alpha_k(t)} 
\partial_{x_k} V_k(t,X^{*}(t)) \right). 
\end{align*}
By definition we have $V_i(T,X)=0$ for all $X.$ Therefore, instead of solving the PMP equation we may ask to solve the $N$ HJB for $V_i=V_i(t,X)$  on 
$[0,T] \times \R^{N}$  for $i=1,\dots,N$ given by the reformulation of the previous equation:
\begin{align}
\label{3:HJB}
\partial_t V_i(t,X) + \sum\limits_{k=1, k\not = i}^{N} \partial_{x_k} V_i(t,X) \left( f_k(X) - \frac{1}{\alpha_k(t)} \partial_{x_k} V_k(t,X) \right) + \partial_{x_i} V_i(t,X) f_i(X) = \\
 - h_i(X) + \frac{1}{2 \alpha_i(t)} (\partial_{x_i} V_i(t,X) )^{2},\nonumber
\end{align}
with terminal condition 
\begin{align}
\label{3:terminal cond}
V_i(T,X) = 0, \; i=1,\dots, N.
\end{align}
\begin{remark}
Since $V_i(t,Y)$ does not contain terminal costs of the type $g_i(X(T))$ the terminal condition 
for $V_i$ is zero. In case of terminal costs we obtain 
$V_i(T,X) = g_i(X(T))$
and terminal constraints for the co-state $\phi^{i}_j$ as 
$\phi^{i}_j(T)= \partial_{x_j} g_i(X^{*}(T))$ in equation \eqref{3:adjoint}.
\par 
The aspect of the game theoretic concept is seen in the HJB equation \eqref{3:HJB} 
in the mixed terms $\partial_{x_k} V_i.$ If we model particles $i$ that   do not
anticipate the optimal choice of the control of other particles $j\not =i,$
then $N$ 
 minimization problems for  $V_i$ in equation \eqref{eq:value fct i} are independent. 
Therefore  the  corresponding HJB for $V_i$ 
and $V_j$ with $j\not =i$ decouple and all mixed terms vanish. In a  different  setting this situation
has been studied in \cite{AlbiHertyPareschi2014aa,AlbiPareschiZanella2014aa} where only a single
control for all particles is present. 
\end{remark}

Assume we have  a (sufficiently regular)
solution $(V_i)_{i=1}^{N}$ with $V_i : [0,T] \times \R^{N} \to \R$. Then, we obtain the optimal control $u_i^{*}(t)$ and the optimal trajectory 
$X^{*}(t)$ for minimizing $V_i$ by 
\begin{equation}\label{3:ctrl thu hjb} 
u_i^{*}(t) = -\frac{1}{\alpha_i(t)} \partial_{x_i} V_i(t,X^{*}(t)), \; i=1,\dots,N,
\end{equation}
where $X^{*}$ fulfills equation \eqref{3:state}. This is an implicit definition of $u_i^{*}.$ However, in view of the dynamics \eqref{3:state} 
it is not necessary to solve equation \eqref{3:ctrl thu hjb} for $u_i^{*}.$ Similar to the discussion in Section \ref{top-left-to-top-right} we obtain 
a controlled systems dynamic provided we have a solution to the HJB equation. The associated controlled dynamics are given by
\begin{equation}\label{3:ctrl HJB}
\frac{d}{dt} x_i(t) = f_i(X(t)) - \frac{1}{\alpha_i(t)} \partial_{x_i} V_i(t,X(t)), \; j=1,\dots,N,
\end{equation}
and initial conditions \eqref{eq:IC}. Comparing the HJB controlled dynamics with equation \eqref{eq:MPC controlled dynamics}
we observe that in the best--reply strategy the full solution to the HJB is not required. Instead, $\partial_{x_i} 
V_i(t,X)$ is approximated by $\partial_{x_i} h_i(X(t)).$  This approximation is also obtained using a discretization 
of equation \eqref{3:HJB} in a MPC framework. Since the equation for $V_i$ is backwards in time we may use
a semi discretization in time on the interval  $(T- \Delta t, T)$ given by 
\begin{align*}
\frac{ V_i(T,X) - V_i(T-\Delta t,X) }{\Delta t} + \sum\limits_{k=1, k\not = i}^{N} \partial_{x_k} V_i(T,X) \left( f_k(X) - \frac{1}{\alpha_k(t)} \partial_{x_k} V_k(T,X) \right) + \\
 + \partial_{x_i} V_i(T,X) f_i(X) =   - h_i(X) + \frac{1}{2 \alpha_i(t)} (\partial_{x_i} V_i(T,X) )^{2} +O(\Delta t), \\
 V_i(T,X)=0.
\end{align*}
Using the terminal condition we obtain that  $V_i(T-\Delta t,X ) = h_i(X)$ for all $X \in \R^N.$ 
\par 
The derivation of the equation for the HJB equation for $V_i(t,Y)$ allows for an arbitrary choice of $T>t.$ 
Hence we may set the terminal time $T$ also to $T:=t+\Delta t.$ This implies to consider the value function 
$$ V_i^{\Delta t}(t,Y) = \int_t^{t+\Delta t} \left( \frac{\alpha_i(s)}2 u_i^{2}(s) + h_i(X(s))  \right) ds$$
where $X(s), s \in (t, t+\Delta t)$ fulfills \eqref{3:state} and where we indicate the dependence on $\Delta t$
by a superscript on $V_i.$ Applying the explicit Euler discretization 
as shown before leads  therefore to 
$$V^{\Delta t}_i(t,Y) =  h_i(Y), \; Y=X(t).$$
Hence, the best--reply strategy applied at time $t$ for a  finite--dimensional problem of $N$ interacting particles
 coincides with an explicit Euler  discretization of the HJB equation for a value function given by $V_i^{\Delta t}(t,Y)$
 where $Y=X(t)$ is the state of the particle system at time $t.$ 

\subsection{ Meanfield limit of the HJB equation \eqref{3:HJB} }
Next, we turn to the meanfield limit of equation \eqref{3:HJB} for $i=1,\dots, N.$ To this end
we assume that  {\bf (A)} and {\bf (B)} holds. We further recall and introduce some notation; 
\begin{align*}
X=(x_i)_{i=1}^{N}, \; Z=(z_i)_{i=1}^{N}, \; \be=(\eta,z_1,\dots,z_{N-1}), \; \be_k:=\left(z_k,  \eta, z_1,\dots, z_{k-1}, z_{k+1}, \dots, z_{N-1} \right).
\end{align*}
\par 
 We obtain the following set of equations for $V_i(t,X)$ and  $i=1,\dots,N,$ 
\begin{align}\label{eq:sym HJB}
\partial_t V_i(t,X) + \sum\limits_{ k=1, k\not = i}^N \partial_{x_k} V_i(t,X) \left( f(x_k,X_{-k}) - \frac{1}{\alpha(t)} \partial_{x_k} V_k(t,X) \right)  \\
+ \partial_{x_i} V_i(t,X) f(x_i, X_{-i}) = - h(x_i,X_{-i}) + \frac{1}{2\alpha(t)} \left( \partial_{x_i} V_i(t,X) \right)^2, \qquad  V_i(t,X)=0.\nonumber
\end{align}
We show that a solution $(V_i)_{i=1}^{N}$ to the previous set of equations is obtained by considering the 
equation \eqref{eq:sol gen} below. 
Suppose that a function $W=W(t,\be): [0,T]\times \R \times \R^{N-1}\to \R$ fulfills 
\begin{align}\label{eq:sol gen}
\partial_t W(t,\be) + \sum\limits_{k=1}^{N-1} \partial_{z_k} W(t,\be)  \left( f( \be_k) - \frac{1}{\alpha(t)} \partial_{\eta} W(t,\be) \right) + \partial_\eta W(t,\be) f(\be) \\
= - h(\be) + \frac{1}{2\alpha(t)}  \left( \partial_\eta W(t,\be) \right)^{2}, \nonumber
\end{align}
and terminal condition $W(T,\be)= 0.$ Suppose a solution $W$ to equation \eqref{eq:sol gen} exists and fulfills the previous equation pointwise a.e. $(t,\be)\in [0,T]\times \R^{N}.$
 Then, we define
\begin{equation}\label{eq:sol spc} 
V_i(t,X) := W(t,x_i,X_{-i}), \; i=1,\dots,N. 
\end{equation}
By definition $W=W(t,\be),$ therefore the partial derivatives 
of $V_i$ are computed as follows where $$(x_i,X_{-i}) = (\eta,z_1,\dots,z_{N-1}):$$ 
\begin{align*}
\partial_t V_i(t,X) &=\partial_t W(t,x_i,X_{-i}), \; \partial_{x_i} V_i(t,X) = \partial_\eta W(t,x_i,X_{-i}), \\
\partial_{x_k} V_k(t,X) &= \partial_{x_k} W(t,x_k,X_{-k}) = \partial_\eta W(t,x_k,X_{-k}), \\
\partial_{x_k} V_i(t,X) &= \partial_{z_{k}} W(t,\be ) \; \mbox{ for } k\in \{1, \dots, i-1 \}, \\
 \partial_{x_k} V_i(t,X) &= \partial_{z_{k-1}} W(t,\be) \; \mbox{ for } k \in \{ i+1, \dots, N\}.   
\end{align*}
 Due to assumption {\bf (A)} we have that 
$$f(\be_k) = f(z_k, z_1,\dots, z_{k-1},z_{k+1}, \dots, z_{i-1}, \eta, z_{i}, \dots, z_{N-1})$$ 
for any $i \in \{1,\dots,N-1 \}.$ The same is true for the argument of $h.$ Therefore, 
$$f(x_k, X_{-k}) = f( \be_k) $$ 
and $h(x_i,X_{-i})=h(\be).$ Therefore, $V_i(t,X)= W(t,x_i,X_{-i})$ fulfills equation
 \eqref{eq:sym HJB}. Hence, instead of studying equation \eqref{eq:sym HJB} we may study the limit for $N\to \infty$
 of equations \eqref{eq:sol gen}. In view of Theorem \ref{Theorem2.1Card} a limit  exists provided
 $W$ is symmetric (and fulfills uniform bound and uniform continuity estimates). 
\par
Note that, $W$ as a solution to equation \eqref{eq:sol gen} is symmetric with respect to the argument $(z_1,\dots,z_{N-1}).$
  This holds true, since $f$ and $h$ are symmetric with respect to $X_{-i}$ for any $i \in \{ 1,\dots, N \}.$  Hence, in the following we may assume to have a solution $W$
to equation \eqref{eq:sol gen} with the property that  for any permutation $\sigma:\{1,\dots,N-1\}\to \{1,\dots,N-1\}$ we have 
\begin{equation}\label{eq:sym of sol}
W(t,\be ) = W(t,\eta, z_{\sigma_1}, \dots, z_{\sigma_{N-1}}).
\end{equation}

 In view of Theorem \ref{Theorem2.1Card} we expect $W(t,\be)$ to converge for for $N\to\infty$
 to a limit function $\W:[0,T] \times \R \times   \mathcal{P}(\R) \to \R$ in the sense of 
 Theorem \ref{Theorem2.1Card}, i.e., up to a subsequence and for $Z \in \R^{N}$ 
\begin{align*}
\lim\limits_{N\to \infty} \sup\limits_{ |\eta|\leq R, t \in [0,T],  Z_{-N} \subset \R^{N-1} } 
| W(t,\be) - {\bf W}(t,\eta,m^{N-1}_{ Z_{-N} } ) | =0
\end{align*} 

Similar to equation \eqref{def:con}  we obtain that the limit 
  $\W:[0,T] \times \R \times \mathcal{P}(\R) \to \R$ fulfills 
the convergence if the measure $m^{N-1}_{Z_{-N}}$ is replaced
by the empirical measure $m^{N}_Z$ for any $Z \in \R^{N}.$    
 Using the introduced notation in Section \ref{sec:setting} we may therefore write 
  $$  W(t,\be)=W^{N}(t,\eta,m^{N-1}_{Z_{-N}}) \sim \W(t,\eta,m^{N}_{ Z } ).$$
\par 
Similarly, we obtain the following meanfield limits for $N$ sufficiently large (and provided the assumptions
of Theorem \ref{Theorem2.1Card} and \cite[Theorem 4.1]{BlanchetCarlier2014aa} are fulfilled. 
\begin{align*}
\partial_t V_i(t,X) &= \partial_t W(t,x_i,X_{-i}) = \partial_t W^{N}(t,x_i,m^{N-1}_{X_{-i}}) & \sim&  \partial_t \W(t,x_i,m^{N}_X), \\
h_i(X) &=h(x_i,X_{-i}) = h^N(x_i,m^{N-1}_{X_{-i}}) &  \sim & \h(x_i,m^{N}_X), \\
(\partial_{x_i} V_i(t,X))^2 &= (\partial_{x_i} W(t,x_i,X_{-i}))^2 = (\partial_{x_i} W^N(t,x_i,m^{N-1}_{X_{-i}}) )^2 & \sim &  (\partial_{x_i} \W(t,x_i,m^{N}_X))^2.
\end{align*}
It remains to discuss the limit of the mixed term 
in equations \eqref{eq:sym HJB} and \eqref{eq:sol gen}, respectively. 
\begin{equation}\label{eq:ttemp} 
\sum\limits_{k=1}^{N-1} \partial_{z_k} W(t,\be) \left( f(\be_k) - \frac{1}{\alpha(t)}  \partial_\eta W(t,\be) \right).
\end{equation}
In order to derive the meanfield limit for equation \eqref{eq:ttemp} we require $f$ to be symmetric
in {\em all} variables, i.e., 
\begin{itemize}
\item[ {\bf (C)} ] We assume $ f(Z) = f( (z_{\sigma_i})_{i=1}^{N} )$ for any permutation $\sigma:\{1,\dots,N\} \to \{1,\dots, N\}$ and for all $Z \in \R^{N}.$  
\end{itemize}
Under assumption {\bf (C)} we have in particular for all $k \in \{1,\dots,N\}$ and a permutation 
$\sigma:\{1,\dots,N-1\} \to \{1,\dots,N-1\}$  
$$ f(\be_k) = f(\eta, z_1,\dots, z_{N-1}) = f(\eta, z_{\sigma_1}, \dots, z_{\sigma_{N-1}} ).$$
Therefore, $f(\be) = f_N(\eta, m^{N-1}_{ Z_{-N} )}).$ In the sense of equation \eqref{def:con}
we further obtain $f_N(\eta,m^{N-1}_{ Z_{-N} )}) \sim {\bf f}(\eta,m^{N}_Z)$ for any $(\eta,Z).$
However under assumption {\bf (C)} we also obtain $f(Z)=f_N(m^{N}_Z) \sim {\bf f}(m^{N}_Z).$ 
Assuming the limit in Theorem \ref{Theorem2.1Card} is unique we obtain that ${\bf f}$ is therefore
{\em independent } of $\eta.$
\par 
 Now, consider the discrete measure $m^{N}_Z=\frac{1}N \sum\limits_{j=1}^N  m^{N}_{z_j}$ and  
 $m_{z_j} = \delta(x-z_j)  \in \mathcal{P}(\R).$ 
For each $j$ we denote by  $m_{z_j}( \zeta) =  \mathcal{Z}(\zeta) \# m_{z_j}$ the push forward
of the discrete measure with the flow field $c:(t,t+a) \times \R  \times \mathcal{P}(\R) \to\R$ and $m_{z_j}(t)=m_{z_j}.$ 
Let the characteristic equations for $\mathcal{Z}$ for fixed $\eta$ be  given by the flow field 
\begin{equation}
\label{eq:psuh c}
\frac{d}{d\zeta} \mathcal{Z}(\zeta) = c(\zeta, \eta, m^{N}_Z(\zeta) ) := \f(  m^{N}_Z(\zeta) )-\frac{1}{\alpha(\zeta)} \partial_\eta \W(\zeta,\eta,m^{N}_Z(\zeta)).
\end{equation}
Similarly to equation \eqref{def:ms},  we obtain the directional derivative of the measure 
of $\W(t,\eta,m^{N}_Z)$with respect to the measure $m^{N}_Z$ in direction of the vectorfield $c$ at $\zeta=t$ 
as 
\begin{align*}
& \sum\limits_{k=1}^{N-1} \partial_{z_k} W(t,\be) \left( f(\be) - \frac{1}{\alpha(t)}  \partial_\eta W(t,\be) \right) \sim 
 \\ 
 \qquad & \langle \partial_m \W(t,\eta,m^{N}_Z), \f(  m^{N}_Z)-\frac{1}{\alpha(t)} \partial_\eta \W(t,\eta,m^{N}_Z) \rangle_{L^2_{m^{N}_Z}},
\end{align*}
where $L^2_{m^{N}_Z}$ denotes the space of square integrable functions for the measure $m^{N}_Z$. 
Performing the limits for $N\to \infty$ in the sense of equation \eqref{def:con}, replacing $\eta$ by $x$, we obtain finally the meanfield
equation for  $\W=\W(t,x,m):[0,T]\times \R \times \mathcal{P}(\R) \to \R$ given by 
\begin{align}\label{eq:meanfield W}
\partial_t \W(t,x,m) + \langle \partial_m \W(t,x,m), \f(m)-\frac{1}{\alpha(t)} \partial_x\W(t,x,m) \rangle_{L^2_m} + \partial_x \W(t,x,m) \F(x,m) \\
 = - \h(x,m) + \frac{1}{2\alpha(t)} \left( \partial_x \W(t,x,m) \right)^{2}, \; 
\W(T,x,m) = 0. \nonumber
\end{align}
The previous equation is reformulated using the concept of directional derivatives of measures $m$ outlined in the Appendix \ref{appendix}. 
Denote by $c(t,x,m) = \f(m) - \frac{1}{\alpha(t)} \partial_x \W(t,x,m)$ a field. 
If $m_{x_j}(t) \in \mathcal{P}(\R)$ for each $t$ is obtained as push forward with the vector field $c,$ then, $m_{x_j}$ fulfills in a weak sense the
continuity equation \eqref{continuity eq}. Therefore, $m^{N}_X= \frac{1}N \sum\limits_{j=1}^{N} m_{x_j}$ fulfills 
\begin{equation}\label{eq:mf t1}
\partial_t m^{N}_X(t,x) + \partial_x \left( c(t,x,m^{N}_X) m^{N}_X(t,x) \right) = 0.
\end{equation}
As seen from the previous equations and the computations in equation \eqref{def:ms}  we therefore have 
 \begin{align*}
  \partial_t \W(t,x,m^{N}_X(t,x)) + \langle \partial_m \W(t,x,m^{N}_X(t,x)), \f(m^{N}_X(t,x)) - \frac{1}{\alpha(t)} \partial_x \W(t,x,m^{N}_X(t,x)) \rangle_{L^2_{m^{N}_X}} = \\
 \frac{d}{dt} \W(t,x,m^{N}_X(t,x)).
  \end{align*}
 This motivates the following definition. For a family of measures $(m(t))_{t \in [0,T]}$ with $m(t,\cdot) \in \mathcal{P}(\R)$, 
define $\v:[0,T]\times \R \to \R$ by 
\begin{equation}
\label{eq:def v}
\v(t,x) := \W(t,x,m(t,x)).
\end{equation}
Then, from equation \eqref{eq:meanfield W} we obtain 
\begin{equation}\label{eq:mf 1}
\partial_t \v(t,x) + \left(  \partial_x \v(t,x)  \right) \f(m) =  - \h(x,m) + \frac{1}{2\alpha(t)} \left( \partial_x \v(t,x) \right)^{2}
 \end{equation}
 and from equation \eqref{eq:mf t1} we obtain using the definition \eqref{eq:def v}
\begin{equation}\label{eq:mf 2}
\partial_t m(t,x) + \partial_x \left(  \left( \f(m) - \frac{1}{\alpha(t)} \partial_x \v(t,x) \right) m(t,x) \right) = 0.
 \end{equation}
 Provided we may solve the meanfield equations \eqref{eq:mf 1} and \eqref{eq:mf 2} for $(\v,m)$ we obtain 
 a solution $\W$ along the characteristics in $m-$space by the implicit relation \eqref{eq:def v}.
In this sense and under the assumptions ${\bf (A)}$ to ${\bf (C)}$ the  meanfield limit of equation \eqref{eq:sym HJB} or respectively equation \eqref{eq:sol gen} is given
  by the system of the following equations \eqref{eq:final mf} and \eqref{eq:final mf2} 
  for $\v:[0,T]\times \R \to \R$ 
and $m(t) \in \mathcal{P}(\R)$ for all $t\in [0,T].$ The terminal condition for $\v$ is given by $\v(T,x)=0.$ 
\begin{align}
\label{eq:final mf}
\partial_t \v(t,x) + \partial_x \left( \v(t,x) \right)  \f(m(t,x)) - \frac{1}{2\alpha(t)} (\partial_x \v(t,x) )^2 &= - \h(x,m(t,x)), \\
\partial_t m(t,x)  + \partial_x \left( \left(\f(m(t,x)) - \frac{1}{\alpha(t)} \partial_x \v(t,x) \right) m(t,x) \right) &= 0. 
\label{eq:final mf2}
\end{align}
We may also express the control $u_i^{*}$ given by equation \eqref{3:ctrl thu hjb}, i.e., 
$$ u_i^{*}(t) = - \frac{1}{\alpha(t)} \partial_{x_i} V_i(t,x_i(t),X_{-i}(t)),$$
in the meanfield limit. Under assumption {\bf (B)} and using equation \eqref{eq:sol spc} and equation \eqref{eq:def v}
For any $X$ we have
\begin{align*}
- \frac{1}{\alpha(t)} \partial_{x_i} V_i(t,X) =  -\frac{1}{\alpha(t)} \partial_x W(t,x_i,X_{-i}) \sim  -\frac{1}{\alpha(t)} \partial_x \W(t,x,m^{N}_X) = -\frac{1}{\alpha(t)} \partial_x \v(t,x).
\end{align*}  

\subsection{ MPC and best reply strategy for the meanfield equation \eqref{eq:final mf}--\eqref{eq:final mf2} }

We obtain the best--reply strategy through a  MPC  approach. Note that the calculations leading to equation \eqref{eq:final mf} 
are independent of the terminal time $T.$ Now, let a time $\tau \in [0,T]$ be fixed
and let $\Delta t>0$ be sufficiently small. Consider the value function on the  receding horizon 
$(\tau, \tau+\Delta t)$ with initial conditions given at $\tau$ and where we, as before,
 add $\Delta t$ as a superscript
to indicate the dependence on the short time horizon:
\begin{equation}\label{eq:vbar}
V^{\Delta t}_i(\tau,Y) = \int_\tau^{\tau+\Delta t} \left( \frac{\alpha_i(s)}2 u_i^{2}(s) + h_i(X(s)) ds \right) ds.
\end{equation} 
 Repeating the derivation of the meanfield limit  computations for $V^{\Delta t}_i$ 
 we obtain  equation \eqref{eq:final mf} defined only for  $t \in [\tau,\tau+\Delta t].$ 
 Also, we obtain  $\v(\tau+\Delta t,x) = 0.$ A first--order in $\Delta t$ approximation of the solution $\v(\tau,x)$ to the (backwards in time )
equation \eqref{eq:final mf} is therefore given by 
\begin{equation}\label{eq:mf app}
\v(\tau,x) =  \h(x,m(t,x)) + O(\Delta t).
\end{equation}
Substituting this relation in the equation for $m$  in \eqref{eq:final mf} 
we obtain the MPC meanfield equation as 
\begin{align}
\label{eq:final mf superfinal}
\partial_t m(t,x)  + \partial_x \left( \left(\f(m(t,x)) - \frac{1}{\alpha(t)} \partial_x \h(x,m(t,x)) \right) m(t,x) \right) &= 0.
\end{align}
This equation is precisely the same as we had obtained for the 
 controlled dynamics using  the best--reply strategy
derived in the previous section and given by equation \eqref{eq:best reply kinetic}. 
\par 

\begin{remark}\label{rem1}
The best--reply strategy for a meanfield game corresponds therefore to considering at each time $\tau$
a value function measuring only the costs for a small next time step. Those costs  may depend
on the optimal choices of the other agents. However, for a small time horizon the derivative of the running costs 
(i.e. $\h$) is a sufficient approximation to the otherwise intractable solution to the full system of meanfield 
equations~\eqref{eq:final mf}-\eqref{eq:final mf2}.
\end{remark}

We summarize the findings in the following Proposition. 
\begin{prop}\label{lemma2}
Assume {\bf (A)} to {\bf (C)} holds true and let $\Delta t>0$ be given. Denote by $\f(m)$ and $\h(x,m)$
the meanfield limit for $N\to\infty$ of $ f(X)$ and $h(X),$ respectively.
Assume that $m:[0,T] \times \R \to \R$ be such that $m(t,\cdot) \in \mathcal{P}(\R)$ and  fulfill equation 
\begin{align}\label{lem:mf 2}
\partial_t m(t,x)  + \partial_x \left( \left(\f(m(t,x)) - \frac{1}{\alpha(t)} \partial_x \h(x,m(t,x)) \right) m(t,x) \right) &= 0.
\end{align}
and let 
\begin{align*}
\v(t,x) = \h(t,x). 
\end{align*}
Then, for any $t\in [0,T]$ and up to an error of order $O(\Delta t)$ the function $\W:[t,t+\Delta t] \times \R \times \mathcal{P}(\R) \to \R$ implicitly defined by 
\begin{align*}
\W(s,x,m(t,x)) = \v(s,x), \; x \in \R, s \in [t,t+\Delta t], 
\end{align*} 
is a solution to the meanfield equation 
\begin{align*}
\partial_s \W(s,x,m) +  \langle \partial_m \W(s,x,m), \f(m)-\frac{1}{\alpha(s)} \partial_x\W(s,x,m) \rangle_{L^2_m} + \partial_x \W(s,x,m) \f(m) \\
 = - \h(x,m) + \frac{1}{2\alpha(s)} \left( \partial_x \W(s,x,m) \right)^{2}, \; \W(t+\Delta t,x,m)=0.
\end{align*}
The meanfield equation is the formal limit for $N\to\infty$ of an $N$ particle game on the time interval $(t,t+\Delta t)$ and described
by equation \eqref{eq:full dynamics} for  $i=1,\dots, N,$ i.e.,
\begin{align*}
\frac{d}{ds} x_i(s) = f_i(X(s)) + u_i(s), \\
u_i(s) = \mbox{ argmin }_{ u:[t,t+\Delta t] \to \R } \int_t^{t+\Delta t} \left( \frac{\alpha_i(s)}2 u^{2}(r) + h_i(X(r)) \right) dr. 
\end{align*}
A solution to the associated $i$th HJB equations for $V_i:[t,t+\Delta t] \times \R^{N} \to \R$ are given by $V_i(t,X):=\W(s,x_i,m^{N}_{X_{-i}})$ for $i=1,\dots,N,$
and the optimal control is $u_i^{*}(s) = -\frac{1}{\alpha_i(s)} \partial_{x_i} V_i(s,X(s)).$ 
\par 
Under assumption {\bf (C)} the meanfield equation \eqref{lem:mf 2} coincides with the 
formal meanfield equation obtained using the best reply strategy \eqref{eq:best reply kinetic}.
\end{prop}

\newpage
\appendix
\section{Technical details}\label{appendix}
We collect some results of \cite{Cardaliaguet2010aa} for convenience. The Kantorowich--Rubenstein 
distance $\done(\mu,\nu)$ for measures $\mu,\nu \in \mathcal{P}(Q)$ is given defined
by 
\begin{equation}
\label{def d1}
\done(\mu,\nu) := \sup \{ \int \phi d(\mu-\nu): \phi:Q\to \R, \phi \mbox{ 1 - Lipschitz } \}. 
\end{equation}

\begin{theorem}[Theorem 2.1\cite{Cardaliaguet2010aa}]\label{Theorem2.1Card}
Let $Q^N$ be a compact subset of $\R^N$. 
 Consider a sequence of functions $(u_N)_{N=1}^\infty$ with $u_N:Q^{N} \to \R.$  
Assume each $u_N(X)=u_N(x_1,\dots,x_N)$ is a symmetric function in all variables, i.e., 
$$ u_N(X)=u_N( x_{\sigma(1)}, \dots, x_{\sigma(N)} )$$
for any permutation $\sigma$ on $\{1,\dots,N\}.$ Denote by
$\done$ the Kantorowich--Rubenstein distance on the space
of probability measures $\mathcal{P}(Q)$ and let $\omega$ 
be a modulus of continuity independent of $N$. Assume that 
 the sequence is uniformly bounded $\| u_N \|_{L^\infty(Q^N)} \leq C$. Further assume that 
  for all $X,Y \in Q^N$ and all $N$  we have
$$ | u_N(X)-u_N(Y) | \leq \omega( {\bf d}_1( m^N_X, m^N_Y) )$$
where $m^N_\xi \in \mathcal{P}(Q)$ is defined by $m^N_\xi(x) = \frac{1}N \sum\limits_{i=1}^N \delta(x-\xi_i)$.
\par 
Then there exists a subsequence $(u_{N_k})_k$ of $(u_N)_N$ and a continuous map 
$U:\mathcal{P}(Q) \to \R$ such that 
\begin{equation}
\label{eq:conv sense}
 \lim\limits_{k\to \infty} \sup\limits_{ X\in \R^N}  | u_{N_k}(X) - U( m^{N_k}_X) | =0.
 \end{equation}
\end{theorem}
An extension is found in \cite[Theorem 4.1]{BlanchetCarlier2014aa}. As toy
example consider $u_N(X)=\frac{1}N \sum\limits_{i=1}^{N} \phi(x_i)$. 
If $\phi:\R\to\R$ is compactly supported, bounded and $|\phi'(\xi)| \leq C$ for all $\xi \in \R,$ 
then the assumptions of the previous theorem are fulfilled. Note that the assumption on $\phi$ implies that for each $i$ we have  $| \partial_{x_i} u_N(X) | \leq \frac{C}N$  for all $X$ and all $N.$
This condition implies the estimate  on $u_N.$  
The corresponding limit is given by the function $U:\mathcal{P}(\R) \to \R$ defined by  
$U(m) = \int \phi dm.$ We have $U(m^{N}_X) = u_N(X).$ 

\par 
Derivatives in the space of measures are described for example in \cite{AmbrosioGigliSavare2008aa}.
They may be motivated by the following formal computation. 
Let $\psi$ be a smooth function on $\R$ and let $y'(t) = c$ for $t\in(a,b)$ and $y(a)=x.$ We denote
by a subindex $t=a$ the evaluation at $t=a$ of the corresponding expression and by a prime
the derivative of $\psi.$
Then,
\begin{align*}
 c \psi'(x) =\int \psi'(z) c \delta( x-z) dz = \left( \int \psi'(z) c \delta(y(t)-z) dz \right)|_{t=a} 
= \\ \left(\int \psi \partial_z \left( c \delta(y(t)-z) \right) \right)|_{t=a}, \\
c \psi'(x) = \left( \frac{d}{dt} \psi(y(t)) \right)|_{t=a} = \left( \frac{d}{dt} \int \psi(z) \delta(y(t)-z) dz \right)|_{t=a}
\end{align*}
Therefore, we may write 
\begin{align*}
\partial_t \delta(y(t)-z) + \partial_z \left( c \; \delta(y(t)-z) \right)=0
\end{align*}
provided that $y'(t)=c$. Further, $\delta(y(t)-z) = y(t) \# \delta(x-z)$ where $\#$ 
is the push forward operator, see below. Hence, for the family of measures $\delta(y(t)-z)$ 
the previous computation lead to a notion of derivatives. This can be formalized to 
a calculus for derivatives in measure space and 
 we summarize in the following more general results from \cite[Chapter II.8]{AmbrosioGigliSavare2008aa}. 
We consider the space of probability measures $\mathcal{P}_p(\R)$ \cite[Equation (5.1.22)]{AmbrosioGigliSavare2008aa}:
$$ \mathcal{P}_p(\R) =\left\{ \mu \in \mathcal{P}(\R): \int |x-\bar{x} |^{p}  d\mu(x) < \infty \mbox{ for some } \bar{x}\in \R \right\}.$$
We assume $\mathcal{P}_p(\R)$ is equipped with
the Wasserstein distance $W_p(\mu,\nu)$ \cite[Chapter 7.1.1]{AmbrosioGigliSavare2008aa}.
In the case $p=1$ and  for bounded measures $\mu,\nu$ this distance is equivalent  to $\done(\nu,\mu)$ defined in equation \eqref{def d1}. For the case $p=2$ we refer to \cite{BenamouBrenier2000aa} for a different characterization. 
\par 
We consider absolutely continuous curves
$m:(a,b)\to \mathcal{P}_p(\R).$ 
The curve $m$ is called absolutely continuous if there exists a function $M\in L^{1}(a,b)$ such that for all  $a\leq s<t \leq b$ we have 
\begin{equation}
\label{abs cont}
W_p( m(s), m(t) ) \leq \int_s^t  M(\xi) d\xi, 
\end{equation}
 see \cite[Definition 1.1.1]{AmbrosioGigliSavare2008aa}.
  For an absolutely continuous curve $m:(a,b)\to\mathcal{P}_p(\R)$, i.e., 
$m(t)\in \mathcal{P}_p(\R),$ and $p>1$ there exists a vector field $v: (a,b) \times 
\mathcal{\R}\to\R$ with $v(t) \in L^p(\R;m(t))$ a.e. $t\in (a,b)$  such that the continuity
equation 
\begin{equation}
\label{continuity eq}
\partial_t m(t,x) + \partial_x \left( v(t,x) m(t,x) \right) = 0
\end{equation}
holds in a distributional sense. Further,  $\| v(t) \|_{L^p(\R;m(t)) } \leq | m' |(t) $ 
a.e. in $t.$ Here, $t\to |m'|(t)$ for $t \in (a,b)$ is the metric derivative of the curve $m.$ 
The precise statement is given in \cite[Theorem 8.3.1]{AmbrosioGigliSavare2008aa} and 
the metric derivative is given in \cite[Theorem 1.1.2]{AmbrosioGigliSavare2008aa} by 
\begin{equation}
\label{metric derv} 
|m'|(t) := \lim\limits_{s\to t} \frac{ W_p( m(s), m(t) ) }{ |s-t| }.
\end{equation}
The limit exists a.e. in $t,$ provided that $m$ is absolutely continuous \eqref{abs cont}. 
 We have $|m'|(t)\leq M(t)$ a.e. for each function $M$ fulfilling 
equation \eqref{abs cont}, see \cite[Chapter 1]{AmbrosioGigliSavare2008aa}.  
Also, the converse result holds true: If $m$ fulfills in a weak sense equation \eqref{continuity eq}
for some $v \in L^1(a,b;L^p(\R;m(\cdot)),$ then $m$ is absolutely continuous. 
Furthermore, solutions to equation \eqref{continuity eq} can be represented 
using the methods of characteristics, see \cite[Lemma 8.1.6, Proposition 8.1.8]{AmbrosioGigliSavare2008aa}. Under suitable assumptions on $m$ and $v$ 
we have that a weak solution to equation \eqref{continuity eq}  is 
\begin{equation}
\label{push forward} m(t,\cdot) = X(t;a,\cdot) \# m(a,\cdot) \; \forall t \in [a,b]
\end{equation}
provided that $X(t)$ solves
 characteristic  system for every $x\in \R$ and every $s \in [a,b]:$
\begin{equation}
\label{characteristic system}
 X(s;s,x) = x \mbox{ and } \partial_t X(t;s,x)= v(t,X(t;s,x)).
\end{equation}
Here, $(s,x)$ is the initial position of the characteristic in phase space and 
$\#$ is the push forward operator, i.e., if applied to 
the set $\{x\}$ we have $m(t,\{X(t;a,x)\}) = m(a,\{X(a;a,x)\})=m(a,\{x\}).$
Equation \eqref{continuity eq} may also be viewed as the directional 
derivative of the family of measures $m(t,\cdot)$ in direction $v.$ 
\par 
In Section \ref{sec:meanfield} we need to discuss a term of the type $\sum\limits_{j=1}^N c(x_j) \partial_{x_j} f(x_1,\dots,x_N)$ for a symmetric function $f.$ Now, consider a family of paths $m_j:(a,b)\to\mathcal{P}(\R)$ generated by  $m_j (t,z) = y_j(t)\#\delta(x_j-z)$
where $y_j$ solves the characteristic equation  $y_j'(t)=c(y_j(t))$ and $y_j(a)=x_j.$ 
Let $m^N_{Y}(t,z) := \frac{1}N \sum\limits_{j=1}^N \delta(y_j(t)-z).$ We have then 
$$\partial_t m^N_Y(t,x) = - \partial_x \left( c(x) m^N_Y(t,x) \right).$$
If we assume that $f$
fulfills the assumption of Theorem \ref{Theorem2.1Card}, then, there exists $\bf:\mathcal{P}(\R) \to \R$
and  then $f(y_1(t),\dots,y_N(t)) = f_{N}( m^{N}_{Y}) \sim  {\bf f}(m^{N}_Y)$ in the sense of
equation \eqref{def:con}. The following computation similar to the motivation shows the 
expression of the unknown term for large $N:$

\begin{align}\label{def:ms}
\sum\limits_{j=1}^{N}  c(x_j) \partial_{x_j} f(x_1,\dots,x_N) = 
\frac{d}{dt} f(y_1(t),\dots, y_j(t), \dots, y_N(t) )|_{t=a}  \\ = \frac{d}{dt} f_{N}(m^{N}_Y(t))|_{t=a} 
\sim \frac{d}{dt} {\bf f}(m^{N}_{Y(t)})|_{t=a}
\end{align}
In order to make the link with the theory developed in \cite{Cardaliaguet2010aa}, we note that the last derivative at $m=m^{N}_Y$ can be interpreted as  
$$\frac{d}{dt} {\bf f}(m^{N}_{Y(t)})|_{t=a} = \langle \partial_m {\bf f}(m), c \rangle_{L^2_m}, $$
with $L^2_m$ the space of square integrable functions with respect to the measure $m$. This formula can either be seen as the definition of $\partial_m {\bf f}(m)$ if one follows the approach of \cite{AmbrosioGigliSavare2008aa} (which is the route taken here) or as a consequence of its definition if one follows the approach of \cite{Cardaliaguet2010aa}. 

\subsubsection*{Acknowledgments}
This work has been supported by KI-Net NSF RNMS grant No. 1107444, grants DFG Cluster of Excellence 'Production technologies for high--wage countries', BMBF KinOpt and CNRS PEPS-HuMaI grant 'DEESSes'. MH and JGL are grateful for the opportunity to stay and work at Imperial College, London, UK, in the Department of Mathematics. PD is on leave 
from CNRS Institut de Math\'{e}matiques de Toulouse UMR 5219; F-31062 Toulouse, France. 
\bibliographystyle{siam}
\bibliography{completeBibTex}
\end{document}